\documentclass[oneside,11pt]{amsart}
\usepackage{amsmath,amsfonts, latexsym,amssymb}
\usepackage[active]{srcltx}
\newtheorem{theorem}{Theorem}[section]

\newtheorem{corollary}[theorem]{Corollary}
\newtheorem{proposition}[theorem]{Proposition}
\newtheorem{conjecture}[theorem]{Conjecture}
\theoremstyle{definition}

\newtheorem{remark}[theorem]{Remark}
\numberwithin{equation}{section}


\begin{document}

\title[A congruence modulo $n^3$ involving two consecutive sums of powers 
\ldots]{A congruence modulo $n^3$ involving two consecutive sums of powers and 
its applications}

\author{Romeo Me\v strovi\' c}

\address{Maritime Faculty, University of Montenegro, Dobrota 36,
 85330 Kotor, Montenegro} \email{romeo@ac.me}

  \begin{abstract}
For various positive integers $k$, the sums of $k$th powers 
of the first $n$ positive integers,
     $$
S_k(n+1)=1^k+2^k+\cdots +n^k,
     $$
have got to be some of the most popular sums in all of mathematics.
In this note we prove that for each $k\ge 2$  
 \begin{equation*}
2S_{2k+1}(n)- (2k+1)nS_{2k}(n)\equiv \left\{
  \begin{array}{ll}
 0\,(\bmod{\,n^3}) & {\rm if}\,\,k\,\,{\rm is\,\,even\,\,or}\,\, n\,\, {\rm 
is\,\, odd}\\
& {\rm or} \,\, n\equiv 0\,(\bmod{\,4})\\
\frac{n^3}{2}\,(\bmod{\,n^3}) & {\rm if}\,\,k\,\,{\rm is\,\, odd} \\
&\,\,{\rm and}\,\, n\equiv 2\,(\bmod{\,4}).
\end{array}\right.
  \end{equation*}

The above congruence allows us to state 
an equivalent formulation  of Giuga's conjecture.
Moreover, we  prove that the first above congruence is satisfied 
modulo $n^4$ whenever  $n\ge 5$ is a prime number such that $n-1\nmid 2k-2$.
In particular, this congruence arises a conjecture for a prime to be 
 Wolstenholme prime. We also propose  several 
Giuga-Agoh's-like conjectures. 
Further, we establish  two congruences modulo 
$n^3$ for two  binomial type sums involving sums of powers $S_{2i}(n)$ 
with $i=0,1,\ldots ,k$. Furthermore, using the above congruence reduced 
modulo $n^2$, we obtain an 
extension of  Carlitz-von Staudt result for odd power sums. 
  \end{abstract}

\maketitle

\vspace{2mm}


\vspace*{1mm}

 \section{Introduction and Basic  Result}

The sum of powers of integers $\sum_{i=1}^ni^k$
 is a well-studied problem in mathematics (see e.g., \cite{bo}, \cite{sc}).
Finding formulas for these sums has interested mathematicians 
for more than 300 years since the time of James Bernoulli (1665-1705).
Several methods were used to find the sum $S_k(n)$ (see for example 
Vakil \cite{va}). These lead to numerous recurrence relations.
For a nice account of  sums of powers see \cite{ed3}.
For simplicity, here as always in the sequel,
 for all integers $k\ge 1$ and $n\ge 2$ we denote
  $$
S_k(n):=\sum_{i=1}^{n-1}i^{k}=1^{k}+2^{k}+3^{k}+\cdots+(n-1)^{k}.
   $$ 
The study of these sums led Jakob Bernoulli
to develop numbers later named in his honor. Namely, the celebrated
{\it Bernoulli's  formula} (sometimes called {\it Faulhaber's formula})
gives the sum $S_k(n)$ explicitly as (see e.g., Beardon \cite{be})
    \begin{equation}
S_k(n)= \frac{1}{k+1} \sum_{i=0}^{k}{k+1\choose i} n^{k+1-i}B_{i}
    \end{equation}
where $B_i$ ($i=0,1,2,\ldots$)
are {\it Bernoulli numbers} defined by the generating function
   $$
\sum_{i=0}^{\infty}B_i\frac{x^i}{i!}=\frac{x}{e^x-1}.
  $$
It is easy to find the values $B_0=1$, $B_1=-\frac{1}{2}$, 
$B_2=\frac{1}{6}$, $B_4=-\frac{1}{30}$, and $B_i=0$ for odd $i\ge 3$. 
Furthermore, $(-1)^{i-1}B_{2i}>0$ for all $i\ge 1$. 
These and many other properties can be found, for instance, in \cite{ir}. 
Several generalizations of the  formula (1.1) were established 
by Z.-H. Sun (\cite[Theorem 2.1]{su1} and \cite{su2}) and 
Z.-W. Sun (\cite{su3}).

By the well known {\it Pascal's identity}
proven by Pascal in 1654 (see e.g., \cite{ms2})
 \begin{equation}
\sum_{i=0}^{k-1}{k\choose i}S_{i}(n+1)=(n+1)^k-1.
   \end{equation}
Recall also that the  formula (1.2) is also presented in Bernoulli's 
{\it Ars Conjectandi} \cite{ber}, 
(also see  \cite[pp. 269--270]{gkp}) 
published posthumously in 1713. 

On the other hand, divisibility properties 
of the sums $S_k(n)$ were investigated by many 
authors (see \cite{ds}, \cite{le}, \cite{ms}, \cite{sm1}). 
For example, in 2003 P. Damianou and P. Schumer 
\cite[Theorem, p. 219]{ds} proved: 

1) if $k$ is odd, then $n$ divides $S_k(n)$ if and only if $n$ is incogruent
to 2 modulo 4; 

2) if $k$ is even, then  $n$ divides $S_k(n)$
if and only if $n$ is not divisible by any prime $p$ such that
$p\mid D_k$, where $D_k$ is the denominator  of the $k$th  Bernoulli
number $B_k$.

Motivated by the recurrence formula for $S_k(n)$
recently obtained by the author in \cite[Corollary 1.7]{me2},
in this note  we prove the following basic result. 
 \begin{theorem}
Let $k$ and $n$ be positive integers. Then for each $k\ge 2$
  \begin{equation}
2S_{2k+1}(n)- (2k+1)nS_{2k}(n)\equiv \left\{
  \begin{array}{ll}
 0\,(\bmod{\,n^3}) & if\,\,k\,\,is\,\,even\,\,or\,\, n\,\, is\,\, odd\\
& or \,\, n\equiv 0\,(\bmod{\,4})\\
\frac{n^3}{2}\,(\bmod{\,n^3}) & if\,\,k\,\,is\,\, odd\\
&and\,\, n\equiv 2\,(\bmod{\,4}).
\end{array}\right.
  \end{equation}
  Furthermore, 
\begin{equation}
2S_{3}(n)-3nS_{2}(n)\equiv\left\{\begin{array}{ll}
 0\pmod{n^3} & if\,\,n\,\, is\,\, odd\\
 \frac{n^3}{2}\pmod{n^3} & if\,\,n\,\, is\,\,even.
\end{array}\right.\end{equation}
In particular, for all $k\ge 1$ and $n\ge 1$, we have
 \begin{equation}
2S_{2k+1}(n)\equiv (2k+1)nS_{2k}(n)\pmod{n^2},
 \end{equation}
and for all $k\ge 1$ and $n\not\equiv 2 (\bmod{\, 4})$
 \begin{equation}
2S_{2k+1}(n)\equiv (2k+1)nS_{2k}(n)\pmod{n^3}.
 \end{equation}
\end{theorem}
 
Combining the congruence (1.5) and the ``even case" 2) of a result of Damianou 
and Schumer \cite[Theorem, p. 219]{ds} mentioned above,
we obtain the following ``odd" extension of their result. 

\begin{corollary}
If $k$ is an odd positive integer and $n$ a positive integer such that 
$n$ is not divisible by any prime $p$ such that
$p\mid D_{k-1}$, where $D_{k-1}$ is the denominator  of the $(k-1)$th  
Bernoulli number $B_{k-1}$. Then  $n^2$ divides $2S_k(n)$.

Conversely, if  $k$ is an odd positive integer and $n$ a positive integer
relatively prime to $k$ such that $n^2$ divides $2S_k(n)$,
then $n$ is not divisible by any prime $p$ such that
$p\mid D_{k-1}$, where $D_{k-1}$ is the denominator  of 
the $(k-1)$th  Bernoulli number $B_{k-1}$.
\end{corollary}

The paper is organized as follows.
Some applications of Theorem 1.1 are presented in the following 
section. In Subsection 2.1 we give  three particular cases 
of the congruence (1.3) of Theorem 1.1 (Corollary 2.1). 
One of these congruences immediately yields a
reformulation  of Giuga's conjecture in terms of the divisibility of 
$2S_n(n)+n^2$ by $n^3$ (Proposition 2.2).

In the next subsection we establish the fact that
the   congruence  (1.6) is satisfied  modulo $n^4$ whenever  
$n\ge 5$ is a prime
number such that $n-1\nmid 2k-2$ ((2.7) of Proposition 2.5). 
Motivated by some particular cases of this congruence and 
related computations via {\tt Mathematica 8}, we propose 
several Giuga-Agoh's-like conjectures. In particular, Conjecture 2.10
characterizes Wolstenholme primes as positive integers $n$ 
such that $S_{n-2}(n)\equiv 0(\bmod{\,n^3})$. 

In Subsection 2.3 we establish  two congruences modulo $n^3$ 
for two  binomial  sums involving sums of powers $S_{2i}(n)$ 
with $i=0,1,\ldots ,k$ (Proposition 2.15).

Combining the congruence (1.5) of Theorem 1.1 with Carlitz-von Staudt result 
for determining $S_{2k}(n)(\bmod{\,n})$ (Theorem 2.18), 
in the last subsection of Section 2 we extend this result modulo $n^2$ for 
power sums $S_{2k+1}(n)$ (Theorem 2.21). We believe that  Theorem 2.21 
can be useful for study some Erd\H{o}s-Moser-like Diophantine equations 
with odd $k$. 

Proofs of all our results are given in Section 3.

\section{Applications of Theorem 1.1}
 \subsection{Variations of Giuga-Agoh's conjecture} 

Taking $k=(n-1)/2$ if $n$ is odd,  $k=n/2$ if $n$ is even
and $k=(n-2)/2$ if $n$ is  even
 into (1.3) of Theorem 1.1,
we respectively obtain the following three congruences.
\begin{corollary} 
If $n$ is an odd positive integer, then 
 \begin{equation}
2S_n(n)\equiv n^2S_{n-1}(n)\pmod{n^3}.
 \end{equation}
If $n$ is even, then 
\begin{equation}
2S_{n+1}(n)- n(n+1)S_n(n)\equiv\left\{\begin{array}{ll}
 0\pmod{n^3} & if\,\,n\equiv 0(\bmod{\,4})\\
 \frac{n^3}{2}\pmod{n^3} & if\,\,n\equiv 2(\bmod{\,4})
\end{array}\right.
  \end{equation}
and 
 \begin{equation}
2S_{n-1}(n)\equiv n(n-1)S_{n-2}(n)\pmod{n^3}.
 \end{equation}
In particular, for each even $n$ we have  
\begin{equation}
S_{n-1}(n)\equiv\left\{\begin{array}{ll}
 0\pmod{n} & if\,\,n\equiv 0(\bmod{\,4})\\
 \frac{n}{2}\pmod{n} & if\,\,n\equiv 2(\bmod{\,4}).
\end{array}\right.
  \end{equation}
\end{corollary}

Notice that if $n$ is any prime, then  by Fermat's little theorem,
$S_{n-1}(n)\equiv -1 (\bmod{\,n})$.  
In 1950 G. Giuga \cite{gi} conjectured that a positive integer
$n\ge 2$ is a prime if and only if 
 $S_{n-1}(n)\equiv -1 (\bmod{\,n})$.
The following proposition provides an equivalent formulation 
of Giuga's conjecture. 

   \begin{proposition}
The following conjectures are equivalent:
 
$(i)$  A positive integer
$n\ge 3$ is a prime if and only if 
 \begin{equation}
S_{n-1}(n)\equiv -1 \pmod{n}.
 \end{equation}

$(ii)$  A positive integer
$n\ge 3$ is a prime if and only if 
 \begin{equation}
2S_n(n)\equiv -n^2\pmod{n^3}.
 \end{equation}
\end{proposition}

Since by  congruence (2.4), $S_{n-1}(n)\not\equiv -1 (\bmod{\,n})$ for each 
even $n\ge 4$, without loss of generality Giuga's conjecture may be 
restricted to the set of odd positive integers.  
In view of this fact and the fact that by (2.1), $n^2\mid S_n(n)$ for each
odd $n$, Proposition 2.2 yields the following equivalent   
formulation of Giuga's conjecture. 
  \begin{conjecture}[Giuga's conjecture]
An odd integer $n\ge 3$ is is a prime if and only if 
 \begin{equation*}
\frac{2S_n(n)}{n^2}\equiv -1\pmod{n}.
 \end{equation*}
 \end{conjecture}

\begin{remark} It is known that $S_{n-1}(n)\equiv -1\, (\bmod{\,n})$
if and only if for each prime divisor $p$ of $n$, $(p-1)\mid (n/p-1)$
and $p\mid (n/p-1)$ (see \cite{gi}, \cite[Theorem 1]{bbb}). Therefore,
any counterexample to Giuga's conjecture must be squarefree. 
Giuga \cite{gi} showed that there are no exceptions to the conjecture
up to $10^{1000}$. In 1985 Bedocchi \cite{bed} improved this bound 
to $n>10^{1700}$. Finally, in 1996 D. Borwein, J.M. Borwein, 
P.B. Borwein and R. Girgensohn raised  the bound to $n>10^{13887}$. 
Recently, F. Luca, C. Pomerance and I. Shparlinski \cite{lps}
proved that for any real number $x$, the number of counterexamples to
Giuga's conjecture $G(x):=\#\{n<x:\, n\,\,\mbox{ is composite and}
\,\,S_{n-1}(n)\equiv -1\, (\bmod{\,n})\}$  
satisfies the estimate $G(x)=o(\sqrt{x})$ as $x\to\infty$.

Independently, in  1990 T. Agoh  
(published in 1995 \cite{ag}; 
see also \cite{bw} and Sloane's sequence A046094 in 
 \cite{sl})  conjectured that a positive integer $n\ge 2$ is a prime if and 
only if  $nB_{n-1}\equiv -1 (\bmod{\,n}).$
Note that the denominator of the  number $nB_{n-1}$ can be greater than 1, but
since by von Staudt-Clausen's theorem (1840) (see e.g., 
\cite[Theorem 118]{hw}; cf. the equality (2.18) given below), the denominator 
of any Bernoulli number $B_{2k}$ is squarefree, it follows that  
the denominator of $nB_{n-1}$ 
is invertible modulo $n$. In 1996  it was  
reported by  T. Agoh \cite{bbb} that his conjecture is equivalent 
to Giuga's conjecture, hence the name Giuga-Agoh's conjecture found in the
litterature.
 It was pointed out in \cite{bbb} 
that this can be seen  from the Bernoulli formula (1.1) after some
analysis involving von Staudt-Clausen's theorem.
The equivalence of both conjectures is in details proved 
in 2002 by B.C. Kellner \cite[Satz 3.1.3, Section 3.1, p. 97]{kel1} 
(also see \cite[Theorem 2.3]{kel2}).  

Quite recently, J.M. Grau and A.M. Oller-Marc\'en 
\cite[Corollary 4]{gm} proved 
that an integer $n$ is a counterexample to Giuga's conjecture
if and only if it is both a Carmichael and a Giuga number
(for definitions and more information on  Carmichael numbers
see W.R. Alford  et al. \cite{agc} and W.D. Banks and C. Pomerance \cite{bp},
and for  Giuga numbers see  D. Borwein et al. \cite{bbb}, 
J.M. Borwein and E. Wong \cite{bw},  and E. Wong \cite[Chapter 2]{w}; 
also see  Sloane's sequences A007850 and A002997  \cite{sl}). 
Furthermore, several open problems concerning Giuga's conjecture
can be found in \cite[8, E Open Problems].
\end{remark}

\subsection{The congruence (1.3) holds modulo $n^4$ for a prime $n\ge 5$}
The following result  shows that for each prime $n\ge 5$ the first congruence
of (1.3) is also satisfied modulo $n^4$.
 \begin{proposition}
Let $p\ge 5$  be a prime  and let $k\ge 2$
be a positive integer such that $p-1\nmid 2k-2$. Then   
      \begin{equation}\label{cong2.4}
2S_{2k+1}(p)\equiv (2k+1)pS_{2k}(p)\pmod{p^4}.
  \end{equation}
Furthermore, if $p-1\nmid 2k$, then    
  \begin{equation}
S_{2k-1}(p)\equiv 0\pmod{p^2}.
  \end{equation}
   \end{proposition}

As a consequence of Proposition 2.5, we obtain the following ``supercongruence"
which generalizes Lemma 2.4 in \cite{me1}.

\begin{corollary} Let $p\ge 5$ 
be a prime  and let $k$ be a positive integer such that 
$k\le (p^4-p^3-4)/2$ and $p-1\nmid 2k+2$.  Then
       \begin{equation}\label{cong2.7}
2R_{2k-1}(p)\equiv (1-2k)pR_{2k}(p)\pmod{p^4}
  \end{equation}
where
  $$
R_{s}(p):=\sum_{i=1}^{p-1}\frac{1}{i^s},\quad s=1,2,\ldots.
  $$
  \end{corollary}
\begin{remark} 
Z.-H. Sun \cite[Section 5, Theorem (5.1)]{su4}
in terms of Bernoulli numbers explicitly determined 
$\sum_{i=1}^{p-1}(1/i^k)\,(\bmod{\, p^3})$ 
for each prime $p\ge 5$ and $k=1,2,\ldots,p-1$. In particular,
substituting the second congruence of Theorem 5.1(a) in \cite{su4}
(with $2k$ instead of even $k$) into (2.7), we immediately obtain 
the following ``supercongruence"
  $$
R_{2k-1}(p)\equiv\frac{k(1-2k)}{2}\left(\frac{B_{2p-2-2k}}{p-1-k}
-4\frac{B_{p-1-2k}}{p-1-2k} \right)p^2\pmod{p^4}
  $$
for all primes $p\ge 7$  and $k=1,\ldots ,(p-5)/2$.

By \cite[(5.1) on page 206]{su4},
    \begin{equation}\label{cong2.10}
S_{2k}(p)\equiv\frac{p}{3}(3B_{2k}+k(2k-1)p^2B_{2k-2})\pmod{p^3},
    \end{equation}
which inserting into (2.7) gives 
   \begin{equation}\label{cong2.9}
S_{2k+1}(p)\equiv\frac{2k+1}{2}p^2B_{2k}\pmod{p^4}
   \end{equation}
for all primes $p\ge 5$ and positive integers $k\ge 2$
such that $p-1\nmid 2k-2$. Moreover, (2.7) with $2k=p-1\ge 4$ 
(i.e., $p\ge 5$) directly gives
   $$
S_{p}(p)\equiv \frac{p^2}{2}S_{p-1}(p)\pmod{p^4}.
  $$
Taking $2k+1=p$ into (2.11), we find that
   $$
S_{p}(p)\equiv\frac{p^3}{2}B_{p-1}\pmod{p^4},
   $$
which  reducing modulo $p^3$, and using the congruence 
$pB_{p-1}\equiv -1(\bmod{\,p})$, yields $2S_{p}(p)\equiv -p^2(\bmod{\,p^3})$.
This is actually  the congruence (2.6) of Proposition 2.2 with a prime 
$n=p\ge 5$. 

Comparing the above two congruences gives 
$S_{p-1}(p)\equiv pB_{p-1}(\bmod{\,p^2})$ for each prime $p\ge 5$.
However, the  congruence (2.10) with $2k=p-1$ implies that 
for all primes $p\ge 5$
  \begin{equation*}
S_{p-1}(p)\equiv pB_{p-1}\pmod{p^3}.
   \end{equation*}
\end{remark}
 \begin{remark}
A computation shows that each of the congruences 
   $$
S_{n}(n)\equiv \frac{n^3}{2}B_{n-1}\pmod{n^4}
   $$
and
  $$
S_{n-1}(n)\equiv nB_{n-1}\pmod{n^3}
   $$
is also satisfied for numerous odd composite positive integers $n$.
However, we propose the following 
 \end{remark}

 \begin{conjecture}
Each of the congruences
$$
S_{n}(n)\equiv \frac{n^3}{2}B_{n-1}\pmod{n^5},
   $$
  $$
S_{n-1}(n)\equiv nB_{n-1}\pmod{n^4}
   $$
 is satisfied for none integer $n\ge 2$.  
\end{conjecture}

Similarly, taking  $k=(p-3)/2$ into (2.11) for each prime $p\ge 5$ 
we get 
    \begin{equation}
S_{p-2}(p)\equiv\frac{(p-2)p^2}{2}B_{p-3}\pmod{p^4}.
   \end{equation}
Therefore, $p^3\mid S_{p-2}(p)$ if and only if the numerator 
of the Bernoulli number $B_{p-3}$ is divisible by $p$,
and such a prime is said to be {\it Wolstenholme prime}.
The only two known such primes  
are 16843 and 2124679, and by a result of McIntosh and 
Roettger from \cite{mr}, these
primes are the only two  Wolstenholme primes less than $10^9$.
In view of the above congruence, and our computation via 
{\tt Mathematica 8} up to $n=20000$ we have the following two conjectures.

 \begin{conjecture}
A positive integer $n\ge 2$ is a Wolstenholme prime if and only if
   $$
S_{n-2}(n)\equiv 0\pmod{n^3}.
   $$
\end{conjecture}

 \begin{conjecture}
The congruence 
   $$
S_{n-2}(n)\equiv 0\pmod{n^4}
   $$
is satisfied for none integer $n\ge 2$. 
\end{conjecture}

\begin{remark} Quite recently, inspired by Giuga's conjecture,  
J.M. Grau, F. Luca and A.M. Oller-Marc\'en 
 \cite{glm} studied the odd positive integers 
$n$ satisfying the congruence 
$$
S_{(n-1)/2}(n)\equiv 0\, (\bmod{\,n}).
 $$
The authors observed \cite[Section 2, Proposition 2.1]{glm} that this 
congruence is satisfied for each odd prime $n$ and for 
each odd positive integer $n\equiv 3\, (\bmod{\,4})$.
Notice that if $n=4k+3$ with $k\ge 0$, the first part of the congruence 
(1.3) yields
  $$
2S_{(n-1)/2}(n)\equiv \frac{(n-1)n}{2}S_{(n-3)/2}(n)\pmod{n^3}
  $$
which is by the congruence (2.7)  satisfied modulo
$n^4$ for each prime $n\ge 7$ such that $n\equiv 3\, (\bmod{\,4})$.
Multiplying the above congruence by $2$ and reducing the modulus,
immediately gives 
 $$
4S_{(n-1)/2}(n)\equiv -nS_{(n-3)/2}(n)\pmod{n^2}.
  $$
The above congruence shows that $S_{(n-1)/2}(n)\equiv 0(\bmod{\,n^2})$  
for some  $n\equiv 3\, (\bmod{\,4})$ if and only if  
$S_{(n-3)/2}(n)\equiv 0(\bmod{\,n})$. 
Furthermore, reducing the congruence (2.11) with $k=(p-3)/4$ 
where $p\ge 7$ 
is a prime such that  $p\equiv 3\, (\bmod{\,4})$ gives
 \begin{equation}
S_{(p-1)/2}(p)\equiv -\frac{p^2}{4}B_{(p-3)/2}\pmod{p^3},
  \end{equation}
whence it follows that for such a prime $p$, 
$S_{(p-1)/2}(p)\equiv 0(\bmod{\,p^2})$.

On the other hand, if $n\equiv 1\, (\bmod{\,4})$, that is
$n=4k+1$ with $k\ge 1$, the first part of the congruence 
(1.3) yields
  $$
2S_{(n+1)/2}(n)\equiv \frac{(n+1)n}{2}S_{(n-1)/2}(n)\pmod{n^3}
  $$
which is by the congruence (2.7) satisfied modulo
$n^4$ for each prime $n\equiv 1\, (\bmod{\,4})$.
Multiplying the above congruence by $2$ and reducing the modulus,
immediately gives 
 $$
4S_{(n+1)/2}(n)\equiv nS_{(n-1)/2}(n)\pmod{n^2}.
  $$ 
The above congruence shows that $S_{(n-1)/2}(n)\equiv 0(\bmod{\,n})$  
for some  $n\equiv 1\, (\bmod{\,4})$ if and only if  
$S_{(n+1)/2}(n)\equiv 0(\bmod{\,n^2})$.  
For example, by \cite[Proposition 2.3]{glm} (cf. Corollary 2.22 given below)
both previous congruences are 
satisfied for every odd prime power $n=p^{2s+1}$ with any  prime 
$p\equiv 1\, (\bmod{\,4})$ and a positive integer $s$. Moreover, 
reducing the congruence (2.10) with $k=(p-1)/4$ where $p\ge 5$ 
is a prime such that  $p\equiv 1\, (\bmod{\,4})$ gives
 \begin{equation}
S_{(p-1)/2}(p)\equiv pB_{(p-1)/2}\pmod{p^2}.
  \end{equation}
The congruence (2.14) shows that 
$S_{(p-1)/2}(p)\equiv 0 (\bmod{\,p^2})$ whenever  $p\equiv 1\, (\bmod{\,4})$ 
is  an irregular prime for which $B_{(p-1)/2}\equiv 0(\bmod{\, p})$.
That the converse is not true shows the composite number 
$n=3737=37\times 101$ satisfying $S_{(n-1)/2}(n)\equiv 0 (\bmod{\,n^2})$
(this is the only such a composite number less than $16000$).   

Nevertheless,  in view of the congruences 
(2.13) and using similar arguments preceeding Conjecture 2.10
(including a computation  up to $n=20000$), 
we have the following  conjecture.
  \end{remark}
 \begin{conjecture}
An odd positive integer $n\ge 3$ such that $n\equiv 3\, (\bmod{\,4})$
satisfies the congruence
   $$
S_{(n-1)/2}(n)\equiv 0\pmod{n^3}
   $$
if and only if $n$ is an irregular prime for which 
$B_{(n-3)/2}\equiv 0(\bmod{\, n})$. 
 \end{conjecture}
We also propose the following 
\begin{conjecture}
The congruence
   $$
S_{(n-1)/2}(n)\equiv 0\pmod{n^3}
   $$
is satisfied for none odd positive integer 
$n\ge 5$ such that $n\equiv 1\, (\bmod{\,4})$.  
   \end{conjecture}

\subsection{Two congruences modulo $n^3$ involving power sums $S_k(n)$}
Combining the congruences of Theorem 1.1 with Pascal's identity,
we can arrive to the following congruences.
 \begin{proposition}
Let $k$ and $n$ be positive integers. Then
 \begin{equation}
2\sum_{i=0}^k(1+n(k+1-i)){2k+2\choose 2i}S_{2i}(n)\equiv -2\pmod{n^3}
  \end{equation}
and 
 \begin{equation}
2\sum_{i=0}^k\left({2k+2\choose 2i}+n(k+1){2k+1\choose 2i}\right)S_{2i}(n)\
\equiv -2\pmod{n^3}.
  \end{equation}
\end{proposition}
 \begin{remark} Clearly, if $n$ is odd, then 
both congruences (2.15) and (2.16) remain also valid 
after dividing  by 2.
However, a computation via {\tt Mathematica 8} for small values $k$
and even $n$ suggests that this would be true  
for each $k$ and even $n$, and so we have
\end{remark}

\begin{conjecture}
The congruence
  \begin{equation*}
\sum_{i=0}^k(1+n(k+1-i)){2k+2\choose 2i}S_{2i}(n)\equiv -1\pmod{n^3}
  \end{equation*}
is satisfied for all $k\ge 1$ and even $n$.
  \end{conjecture}
\subsection{An extension of Carlitz-von Staudt result for odd power sums}
The following congruence is known as a {\it Carlitz-von Staudt's result}
 \cite{ca} in  1961  (for an easier proof see \cite[Theorem 3]{mo3}).
  \begin{theorem}{\rm(\cite{ca}, \cite[Theorem 3]{mo3})} 
Let $k$ and $n>1$ be positive integers. Then
   \begin{equation}
S_k(n)\equiv\left\{
  \begin{array}{ll}
0 \pmod{\frac{(n-1)n}{2}} & if\,\, k\,\, is\,\, odd\\
-\sum_{(p-1)\mid k, p\mid n}\frac{n}{p}\pmod{n} & 
if\,\, k\,\, is\,\, even
\end{array}\right.\end{equation}
where the summation is taken over all primes $p$ such that 
$(p-1)\mid k$ and $p\mid n$.
  \end{theorem}
     \begin{remark}
It is easy to show the first (``odd") part of Theorem 2.18 
(see e.g., \cite[Proof of Theorem 3]{mo3} or
  \cite[Proposition 1]{ms}) whose proof is a modification of 
Lengyel's arguments in \cite{le}. 
Recall also that the classical {\it theorem of Faulhaber} states that every 
sum $S_{2k-1}(n)$ (of odd power) can be expressed as a polynomial of the 
triangular number $T_{n-1}:=(n-1)n/2$; see e.g.,  \cite{be} 
or \cite{kn}.
For even powers, it has been shown  that the sum $S_{2k}(n)$
is a polynomial in the triangular number $T_{n-1}$ multiplied 
by a linear factor in $n$ (see e.g., \cite{kn}).
\end{remark}
  \begin{remark}
The second part of the congruence (2.17) in Theorem 2.18 can be proved using  
the famous von Staudt-Clausen theorem (given below) 
discovered  without proof
by C. Clausen \cite{cl} in 1840, and independently by K.G.C. von Staudt 
in 1840 \cite{st}; for alternative proofs, see, e.g., Carlitz \cite{ca}, 
Moree \cite{mo1} or Moree \cite[Theorem 3]{mo3}.
This also follows from Chowla's proof of von Staudt-Clausen theorem 
\cite{ch}.
Namely, Chowla proved that the difference 
    $$
\frac{S_{2k}(n+1)}{n}-B_{2k}
    $$
is an integer for all positive integers $k$ and $n$.  This 
together with the facts  that $S_{2k}(n+1)\equiv S_{2k}(n)\,(\bmod{\, n})$
and that by  {\it von Staudt-Clausen theorem},
     \begin{equation}
B_{2k}=A_{2k}-\sum_{(p-1)\mid 2k\atop p\,\,\mathrm{prime}}\frac{1}{p},
    \end{equation}
where $A_{2k}$ is an integer, immediately gives the second part of 
the congruence (2.17).

Recall also that in many places,  
the von Staudt-Clausen theorem is stated in the following equivalent 
statement (e.g., see \cite[page 153]{su}):
  \begin{equation*}
pB_{2k}\equiv \left\{
  \begin{array}{ll}
\,\,\,\, 0\,(\bmod{\,p}) & {\rm if}\,\,p-1\nmid 2k\\
-1\,(\bmod{\,p}) & {\rm if}\,\,p-1\mid 2k,
\end{array}\right.
  \end{equation*}
where $p$ is a prime and $k$ a positive integer.
      \end{remark}

Combining the congruence (1.5) in Theorem 1.1  
with the second (``even") part of the congruence (2.17),  
 we immediately obtain an improvement  of its first (``odd")  part as follows.
  \begin{theorem}
Let Let $k$  and $n$ be positive integers. Then
    \begin{equation}
2S_{2k+1}(n)\equiv -(2k+1)n\sum_{(p-1)\mid 2k, p\mid n}\frac{n}{p}\pmod{n^2},
   \end{equation}
where the summation is taken over all primes $p$ such that 
$p-1\mid k$ and $p\mid n$.
   \end{theorem}
In particular, taking $n=p^s$ and $k=(p^s-1)/4$ into (2.19)
where $p$ is an odd prime $p$ and $s\ge 1$ such that 
$p^s\equiv 1(\bmod{\, 4})$, we immediately obtain an 
  analogue of Proposition 2.3 in a recent paper \cite{glm}. 
  \begin{corollary}
Let $p$ be an odd prime. Then 
 \begin{equation*}
S_{(p^s+1)/2}(p^s)\equiv \left\{
  \begin{array}{ll}
\,\,\, 0\,(\bmod{\,p^{2s}}) & if \,\,
p\equiv 1(\bmod{\, 4})\,\, and\,\, s\ge 1\,\,\, is\,\,\, odd\\
 -\frac{p^{2s-1}}{4}\,(\bmod{\,p^{2s}}) &  if\,\,\, s\ge 2\,\,\, is\,\,\,
even.\\
 \end{array}\right.
    \end{equation*}
  \end{corollary}

Finally, comparing (2.17), (2.18) and (2.19), we  immediately obtain 
an ``odd" extension of a result due to Kellner \cite[Theorem 1.2]{kel2}
in 2004 (the congruence (2.20) given below).
   \begin{corollary}
Let Let $k$  and $n$ be positive integers. Then
    \begin{equation}
S_{2k}(n)\equiv nB_{2k}\pmod{n} \,\,\,\,(Kellner \,\, [20])
  \end{equation}
and
      \begin{equation}
2S_{2k+1}(n)\equiv (2k+1)n^2B_{2k}\pmod{n^2}.
   \end{equation}
\end{corollary}

\begin{remark}
Notice also that Theorem 2.18  plays a key role in a recent  study 
(\cite{mo1}, \cite{mo3}, \cite{gmz}, \cite{mru})
of the {\it Erd\H{o}s-Moser Diophantine equation} 
      \begin{equation}
1^k+2^k+\cdots +(m-2)^k+(m-1)^k=m^k
     \end{equation}
where $m\ge 2$ and $k\ge 2$ are positive integers.
Notice that $(m,k)=(3,1)$ is the only solution for $k=1$.
In letter to Leo Moser around 1950, Paul Erd\H{o}s
conjectured that such solutions of the above equation do not exist
(see \cite{mos}). Using remarkably elementary 
methods, Moser \cite{mos} showed in 1953 that if 
$(m,k)$ is a solution of (2.22) with $m\ge 2$ and $k\ge 2$, 
then $m>10^{10^6}$ and $k$ is even. 
 Recently,  using Theorem 2.18,  P. Moree 
\cite[Theorem 4]{mo3} improved  the bound 
on $m$ to $1.485\cdot 10^{9321155}$.
  That Theorem 2.18 can be used 
to reprove Moser's result was first observed in 1996 by Moree \cite{mo2}, 
where it played a key role in the study of the more general 
equation 
    \begin{equation}
1^k+2^k+\cdots +(m-2)^k+(m-1)^k=am^k
  \end{equation}
where $a$ is a given positive integer. 
Moree \cite{mo2} generalized Erd\H{os}-Moser 
conjecture in the sense that the only solution 
of the ``generalized" Erd\H{os}-Moser Diophantine
equation (2.23) is the trivial solution $1+2+\cdots +2a=a(2a+1)$.
Notice also that Moree \cite[Proposition 9]{mo2} proved that in any 
solution of the equation (2.23), $m$ is odd.      
Nevertheless, motivated by the Moser's technique  
\cite[proof of Theorem 3]{mo3}  previously  mentioned,
to study  (2.22),  we believe that  Theorem 2.21 
can be useful for study  some other Erd\H{o}s-Moser type Diophantine equations 
with odd $k$. 
  \end{remark}

\section{Proofs of Theorem 1.1, Corollaries 1.2, 2.6
and Propositions 2.2, 2.5 and 2.15} 

\begin{proof}[Proof of Theorem $1.1$]
If $k\ge 1$ then by the binomial formula, for each $i=1,2,\ldots ,n-1$
we have  
  \begin{equation}\begin{split}
& 2(i^{2k+1}+(n-i)^{2k+1})-(2k+1)n(i^{2k}+(n-i)^{2k})\\
\equiv & 2\left(i^{2k+1}-i^{2k+1})+{2k+1\choose 1}ni^{2k}
-{2k+1\choose 2}n^2i^{2k-1}\right)\\
&-(2k+1)n\left(i^{2k}+i^{2k}-{2k\choose 1}ni^{2k-1}\right)\pmod{n^3}\\
=& 2(2k+1)ni^{2k}-2(2k+1)kn^2i^{2k-1}-2(2k+1)ni^{2k}+
2(2k+1)kn^2i^{2k-1}\\
=&0\pmod{n^3}.
   \end{split}\end{equation}
If $k\ge 3$ and  $n$ is odd then after summation of (3.1) 
over $i=1,2,\ldots,(n-1)/2$ we obtain
   \begin{equation}
2\sum_{i=1}^{n-1}i^{2k+1}-(2k+1)n\sum_{i=1}^{n-1}i^{2k}\equiv 0\pmod{n^3}.
   \end{equation}
If $k\ge 2$ and $n$ is even
then after summation of (3.1) over $i=1,2,\ldots,n/2$ we get
   \begin{equation*}
2\sum_{i=1}^{n-1}i^{2k+1}+
2\left(\frac{n}{2}\right)^{2k+1}-
(2k+1)n\sum_{i=1}^{n-1}i^{2k}- 
(2k+1)n\left(\frac{n}{2}\right)^{2k}\equiv 0\pmod{n^3},
   \end{equation*}
or equivalently,
   \begin{equation}
2S_{2k+1}-(2k+1)nS_{2k}\equiv \frac{kn^{2k+1}}{2^{2k-1}}=
\frac{n^3}{2}\cdot k\cdot \left(\frac{n}{2}\right)^{2k-2}
\pmod{n^3}.
   \end{equation}
Since for even $n$
 \begin{equation*}
\frac{n^3}{2}\cdot k\cdot \left(\frac{n}{2}\right)^{2k-2}\equiv \left\{
  \begin{array}{ll}
 0\,(\bmod{\,n^3}) & {\rm if}\,\,k\,\,{\rm is}\,\,{\rm even\,\,or} 
\,\, n\equiv 0\,
(\bmod{\,4})\\
\frac{n^3}{2}\,(\bmod{\,n^3}) & {\rm if}\,\,k\,\,{\rm is\,\, odd \,\,and}
\,\, n\equiv 2\,(\bmod{\,4}),
\end{array}\right.
  \end{equation*}
this together with (3.3) and (3.2) yields both congruences of (1.3) in 
Theorem 1.1.

Finally, for $k=1$ we have
  \begin{equation*}
2S_{3}(n)-3nS_{2}(n)= \frac{n^3}{2}\cdot 
(1-n)\equiv\left\{\begin{array}{ll}
 0\pmod{n^3} & {\rm if}\,\,n\,\, {\rm is\,\, odd}\\
 \frac{n^3}{2}\pmod{n^3} & {\rm if}\,\,n\,\, {\rm is\,\,even}.
\end{array}\right.
   \end{equation*}
This completes the proof.
   \end{proof}
\begin{proof}[Proof of Corollary $1.2$] 
Both assertions immediately follow applying the congruence (1.5) 
and a result of P. Damianou and P. Schumer 
  \cite[Theorem, p. 219]{ds} 
which asserts that if $k$ is even, then  $n$ divides $S_k(n)$
if and only if $n$ is not divisible by any prime $p$ such that
$p\mid D_k$, where $D_k$ is the denominator  of the $k$th  Bernoulli
number $B_k$.
 \end{proof}
\begin{proof}[Proof of Proposition $2.2$] {\it Proof of} 
$(i)\Rightarrow (ii)$. Suppose that Giuga's conjecture is true. Then if $n$ 
is an odd positive integer satisfying the congruence (2.6) of Proposition 2.2,
using this and (2.1) of Corollary 2.1, we find that   
  $$
n^2S_{n-1}(n)\equiv 2S_n(n)\equiv -n^2\pmod{n^3},
  $$
whence we have 
  $$
S_{n-1}(n)\equiv -1\pmod{n}.
  $$
By Giuga's conjecture, the above congruence implies that
$n$ is a prime.

If  $n\ge 4$ is an even positive integer, then the congruence (2.4) shows 
that $S_{n-1}(n)\not\equiv -1 (\bmod{\,n})$. 
We will show that for such a $n$, $2S_{n}(n)\not\equiv -n^2(\bmod{\,n^3})$.
Take $n=2^s(2l-1)$, where $s$ and $l$ are positive integers. Since 
for $i=1,2,\ldots$ we have $(2i)^n\equiv 0 (\bmod{\,2^n})$,
this together with the inequality $2^{2^s}\ge 2^{s+1}$ yields 
$(2i)^n\equiv 0 (\bmod{\,2^{s+1}})$. Therefore, we obtain
 $$
2S_n(n)\equiv 2\sum_{1\le j\le n-1\atop j\, {\rm odd}}j^n\pmod{2^{s+1}}. 
  $$ 
Since by Euler theorem, for each odd $j$ 
  $$
j^n=j^{2^s(2l-1)}=\left(j^{2^s}\right)^{2l-1}=
\left(j^{\varphi(2^{s+1})}\right)^{2l-1}\equiv 1\pmod{2^{s+1}}
\equiv 1\pmod{2^s},
  $$
where $\varphi(m)$ is the Euler's totient function, 
then substituting this into above congruence, we get
  $$
2S_n(n)\equiv n=2^s(2l-1)\not\equiv 0\pmod{2^{s+1}}. 
  $$ 
Now, if  we suppose that $2S_{n}(n)\equiv -n^2(\bmod{\,n^3})$,
then must be $2S_{n}(n)\equiv 0(\bmod{\,n^2})$, and so,
$2S_{n}(n)\equiv 0(\bmod{\,2^{2s}})\equiv 0
(\bmod{\,2^{s+1}})$. This contradicts the above congruenece,
and the impication $(i)\Rightarrow (ii)$ is proved.

{\it Proof} of $(ii)\Rightarrow (i)$. 
Now suppose that Conjecture (ii) of Proposition 2.2 is true. Then if $n$ 
is an odd positive integer satisfying the congruence (2.5), multiplying this 
by $n^2$ and using (2.1) of Corollary 2.1, we find that   
  $$
2S_n(n)\equiv n^2S_{n-1}(n)\equiv -n^2\pmod{n^3},
  $$
which implies that 
  $$
2S_n(n)\equiv -n^2\pmod{n^3}.
  $$
By our Conjecture (ii), the above congruence implies that
$n$ is a prime.

If  $n\ge 4$ is an even positive integer, then we have previously shown that 
for such a $n$, $2S_{n}(n)\not\equiv -n^2(\bmod{\,n^3})$ 
and $S_{n-1}(n)\not\equiv -1 (\bmod{\,n})$. 
This completes the proof of  impication $(ii)\Rightarrow (i)$.
  \end{proof}

\begin{proof}[Proof of Proposition $2.5$]
If we extend the congruence (3.1) modulo $n^4$,  
then in the same manner we obtain 
 \begin{equation*}\begin{split}
& 2(i^{2k+1}+(n-i)^{2k+1})-(2k+1)n(i^{2k}+(n-i)^{2k})\\
&\equiv 2{2k+1\choose 3}n^3i^{2k-2}-
(2k+1){2k\choose 2}n^3i^{2k-2}\pmod{n^4},
  \end{split}\end{equation*}
whence it follows that
\begin{equation}
2S_{2k+1}(n)-(2k+1)nS_{2k}(n)\equiv \frac{k(1-4k^2)}{3}n^3S_{2k-2}(n)
\pmod{n^4}.
  \end{equation}
If $n=p$ is a prime such that $p-1\nmid 2k-2$,
then the well  known congruence $S_{2k-2}(p)\equiv 0\,(\bmod{\,p})$
(see e.g., \cite[the congruence (6.3)]{su1} or
 \cite[Theorem 1]{ms2}) and (3.4) yield the congruence (2.7).
 Finally, (2.8) immediately follows reducing  (2.7) modulo
$p^2$ and using the previous fact that  $S_{2k}(p)\equiv 0\,(\bmod{\,p})$
whenever $p-1\nmid 2k$.
 \end{proof}
\begin{remark}
Applying a  result of P. Damianou and P. Schumer \cite[Theorem, p. 219]{ds} 
used in the proof of Proposition 2.5 to the congruence (3.4),
it follows that
      $$
2S_{2k+1}(n)\equiv  (2k+1)nS_{2k}(n)\pmod{n^4}
      $$
whenever $n$ is not divisible by any prime $p$ such that
$p\mid D_{2k-2}$, where $D_{2k-2}$ is the denominator  of the $(2k-2)$th  
Bernoulli number $B_{2k-2}$. The converse assertion is true if $n$ is 
relatively prime to the integer $k(1-4k^2)/3$. 
\end{remark}
 \begin{proof}[Proof of Corollary $2.6$] 
By Euler's theorem \cite{hw}, for all positive integers $m$ and $i$
such that  $1\le m< p^4-p^3$ and $1\le i\le p-1$ we have 
   $1/i^m\equiv i^{\varphi(p^4)-m}\,\,(\bmod{\,\,p^4})$,
where $\varphi(p^4)=p^4-p^3$ is the Euler's totient function.
Therefore, $R_m\equiv S_{p^4-p^3-m}\,(\bmod{\,p^4})$.
Applying the last congruence for $m=2k-1$ and $m=2k$,
and substituting this into (2.7)  of Proposition 2.5
 with $p^4-p^3-2k\ge 4$ instead of $2k$, we  immediately obtain 
  $$
2R_{2k-1}(p)\equiv (p^4-p^3-2k+1)pR_{2k}(p)\equiv
(1-2k)pR_{2k}(p)\pmod{p^{4}},
 $$
as desired.
  \end{proof}
\begin{proof}[Proof of Proposition $2.15$]  As 
$S_0(n)=n-1$ and $S_1(n)=(n-1)n/2$, Pascal's identity (1.2) yields
 \begin{equation}\begin{split}
&2(n^{2k+2}-1) =2\sum_{i=0}^{2k+1}{2k+2\choose i}S_{i}(n)\\
&=2(n-1)(1+(k+1)n)+\sum_{i=1}^{k}\left(2{2k+2\choose 2i}S_{2i}(n)+
 2{2k+2\choose 2i+1}S_{2i+1}(n) \right).
   \end{split}\end{equation}
If $n$ is odd, then multiplying the congruence 
(1.6) of Theorem 1.1 by ${2k+2\choose 2i+1}$
and using  the identity ${2k+2\choose 2i+1}=
\frac{2k+2-2i}{2i+1}{2k+2\choose 2i}$, we find that
 \begin{equation}\begin{split}
{2k+2\choose 2i+1}2S_{2i+1}(n)
&\equiv \frac{2k+2-2i}{2i+1}{2k+2\choose 2i}(2i+1)nS_{2i}(n)\pmod{n^3}\\
&= (2k+2-2i){2k+2\choose 2i}nS_{2i}(n)\pmod{n^3}
   \end{split}\end{equation}
for each $i=1,\ldots,k$.
Now substituting (3.6) into (3.5), we obtain
\begin{equation}
2(n-1)(1+(k+1)n)+2\sum_{i=1}^k(1+n(k+1-i)){2k+2\choose 2i}S_{2i}(n)
\equiv -2\pmod{n^3},
  \end{equation}
which is obviously the same as (2.15).

If $n$ is even, then since  ${2k+2\choose 2i+1}$ is even (this is true 
by the identity ${2k+2\choose 2i+1}=\frac{2(k+1)}{2i+1}{2k+1\choose 2i}$),
we have   that ${2k+2\choose 2i+1}\frac{n^3}{2}\equiv 0\,(\bmod\,{n^3})$.
This shows that (2.15) is satisfied for  even $n$ and each $i=1,\ldots,k$, 
and hence, proceeding in the same manner as in the previous case,
we obtain (2.15). 

Further, applying the identities $2i{2k+2\choose 2i}=(2k+2){2k+1\choose 2i-1}$
and ${2k+2\choose 2i}-{2k+1\choose 2i-1}={2k+1\choose 2i}$,
the left hand side of (2.16)  is equal to
\begin{equation*}\begin{split}
&2(1+n(k+1))\sum_{i=0}^k{2k+2\choose 2i}S_{2i}(n)-
n\sum_{i=0}^k2i{2k+2\choose 2i}S_{2i}(n)\\
=&2\sum_{i=0}^k{2k+2\choose 2i}S_{2i}(n)+2n(k+1)(n-1)+
2n(k+1)\sum_{i=1}^k{2k+2\choose 2i}S_{2i}(n)\\
&-2n(k+1)\sum_{i=1}^k{2k+1\choose 2i-1}S_{2i}(n)\\
=&2\sum_{i=0}^k{2k+2\choose 2i}S_{2i}(n)+2n(k+1)(n-1)\\
&+2n(k+1)\sum_{i=1}^k\left({2k+2\choose 2i}- {2k+1\choose 2i-1}
\right)S_{2i}(n)\\
=&2\sum_{i=0}^k{2k+2\choose 2i}S_{2i}(n)+2n(k+1)(n-1)+
2n(k+1)\sum_{i=1}^k{2k+1\choose 2i}S_{2i}(n)\\
=&2\sum_{i=0}^k{2k+2\choose 2i}S_{2i}(n)+2n(k+1)
\sum_{i=0}^k{2k+1\choose 2i}S_{2i}(n)\\
=&2\sum_{i=0}^k\left({2k+2\choose 2i}+n(k+1){2k+1\choose 2i}\right)S_{2i}(n).
  \end{split}
 \end{equation*}
Comparing the above equality with (2.15) immediately gives (2.16).
  \end{proof}

\end{document}